\documentclass[12pt,english]{amsart}
\usepackage{pslatex}
\usepackage[T1]{fontenc}
\usepackage[latin1]{inputenc}
\usepackage{graphicx}
\usepackage{setspace}
\usepackage{amssymb}
\usepackage[numbers]{natbib}

\makeatletter
 \theoremstyle{plain}    
 \newtheorem{thm}{Theorem}[section]
 \numberwithin{equation}{section} 
 \numberwithin{figure}{section} 
 \theoremstyle{plain}
 \theoremstyle{remark}
 \newtheorem{rem}[thm]{Remark}
 \theoremstyle{plain}    
 \newtheorem{prop}[thm]{Proposition} 
 \theoremstyle{plain}    
 \newtheorem{cor}[thm]{Corollary} 

\usepackage{babel}
\makeatother
\begin{document}

\title{Epidemic branching processes with and without vaccination}

\author{\textbf{Arni S.R. Srinivasa Rao$^{*}$ and Chris T. Bauch}}

\maketitle
Department of Mathematics and Statistics, University of Guelph, Canada,
N1G 2W1. 

{*} Present address: Mathematical Institute, 24-29 St Giles, University
of Oxford, Oxford, OX1 3LB England. Emails: <arni@maths.ox.ac.uk>,
<cbauch@uoguelph.ca>

\begin{abstract}
Here we treat the transmission of disease through a population as
a standard Galton-Watson branching process, modified to take the presence
of vaccination into account. Vaccination reduces the number of secondary
infections produced per infected individual. We show that introducing
vaccination in a population therefore reduces the expected time to
extinction of the infection. We also prove results relating the distribution
of number of secondery infections with and without vaccinations. 
\end{abstract}



\section{Introduction and background}

After the introduction of penicillin and mass vaccination campaigns
of the 1950s, it was thought that infectious diseases would soon be
a thing of the past. However, infectious diseases have turned out
to be an intractable problem, as new infectious diseases continually
emerge while previously existing diseases circumvent existing control
measures through evolution. Accordingly, there is a continuing need
to develop mathematically rigorously methods that can be applied to
infectious disease epidemiology, in order to predict future spread
of infections and assess the likely impacts of alternative control
measures, including vaccination.

One such area of promise emerging from the field of probability is
branching processes. Infectious disease transmission can be thought
of as a branching process in which each infected individual in a given
generation of infection produces some number of infected individuals
in the next generation, as described by a probability density function.
This picture is particularly apt in the early stages of an outbreak,
when the proportion of infected individuals is small and the process
is highly stochastic. Epidemiologists need to know whether a disease
will die out, and how many individuals are likely to be infected before
the extinction, given certain fundamental probability distributions
such as the distribution of number of secondary infections produced
per infected individual. They also need to know how these parameters
will change once control measures such as vaccination are introduced.

Some work has started in applying branching process-type arguments
successfully to the mathematical modelling of disease transmission
\cite{key-0}. However much work remains to be done in developing
a mathematically rigorous methodology to understand impact of vaccines.
Branching processes have been applied to areas in biology such as
the study of offspring distributions and the extinction of family
surnames. However, there are issues that are particular to infectious
disease epidemiology, such as the impact of vaccination on the branching
processes. These pecularities necessitate an extension of branching
process methodologies to the case of infectious disease transmission
and vaccination. Here, infectious disease transmission and vaccination
are treated as a classic Galton-Watson branching process. 

Let $Z_{n}$ be the number of infected individuals in the $n^{th}$
generation of disease transmission, where $n\geq0$. Since we consider
the early stages of disease transmission and the proportion of susceptibles
is almost 1, infected individuals do not co-infect the same susceptible
individuals. An `infected' (resp. `susceptible') here means an individual
who is infected (resp. susceptible) in the present generation. Every
infected individual has the same probability distribution of the number
of individuals they infect, and these are distributed independently
of one other. Assume each infected in the $n^{th}$ generation generates
$N$ new individuals in the $(n+1)^{th}$ generation. Here $N\in\mathbb{Z^{+}}$is
a measurable function (a r.v.) whose infected population distribution
is $\left\{ \theta_{i}\right\} _{i=0}^{\infty}$ i.e. $\Pr\left[N=i\right]=\theta_{i},\, i\geq0.$

Let $N_{r}^{(n+1)}$ represents the number of infected by the $r^{th}$
individual in $n^{th}$ generation. These infected are account for
the size of the population in the $(n+1)^{th}$ generation. Then $Z_{n+1}=N_{1}^{(n+1)}+N_{2}^{(n+1)}+...+N_{Z_{n}}^{(n+1)}.$
Assume that each $N_{r}^{(n)}$ is \emph{iid} r.v.'s i.e. $\Pr\left[N_{r}^{(n)}=i\right]=\Pr\left[N=i\right]=\theta_{i}.$
The size of the $(n+1)^{th}$ generation $Z_{n+1}$ depends on the
number of infected individuals in the $n^{th}$ generation, therefore
$Z_{n}$ satisfies the properties of a Markov chain. We assume that
$Z_{0}=1$, guaranteeing that there is one infected individual from
which the epidemic can start. After one disease generation (one discrete
time interval), this individual has infected a further number of individuals
as given by the probability distribution$\{\theta_{i}\}_{i=0}^{\infty}$.
By this time, the original infected individual has also recovered
(or died) from the infection. Thus the first generation constitutes
$Z_{1}$ infecteds, which is the sum of the $Z_{0}$ random variables
each with a probability $\theta_{i}$. The second generation constitutes
$Z_{2}$ infected and the process continues. We call this process
an epidemic branching process (EBP).

\begin{figure}
\includegraphics[%
  scale=0.4]{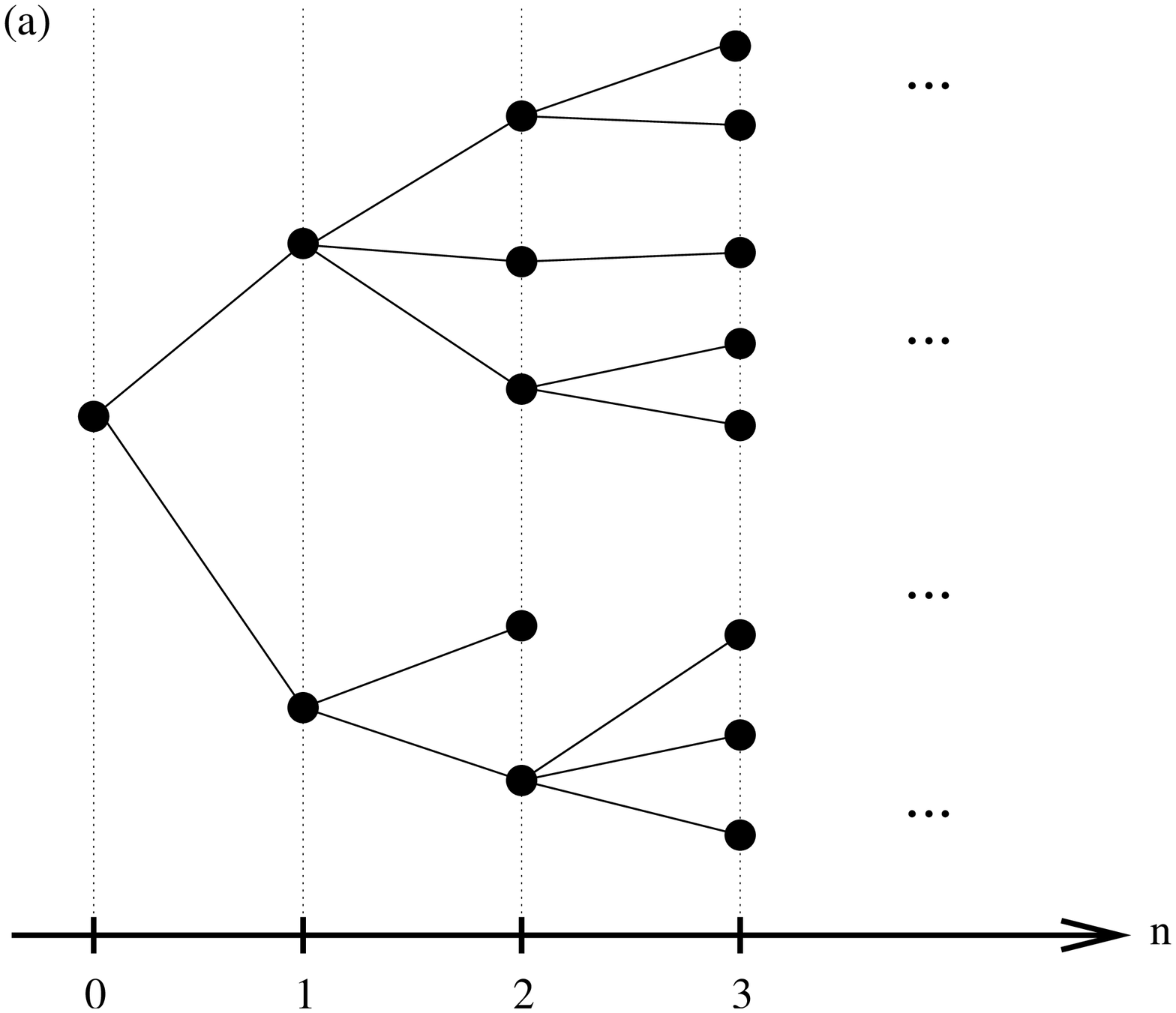}

$ $

$ $

\includegraphics[%
  scale=0.4]{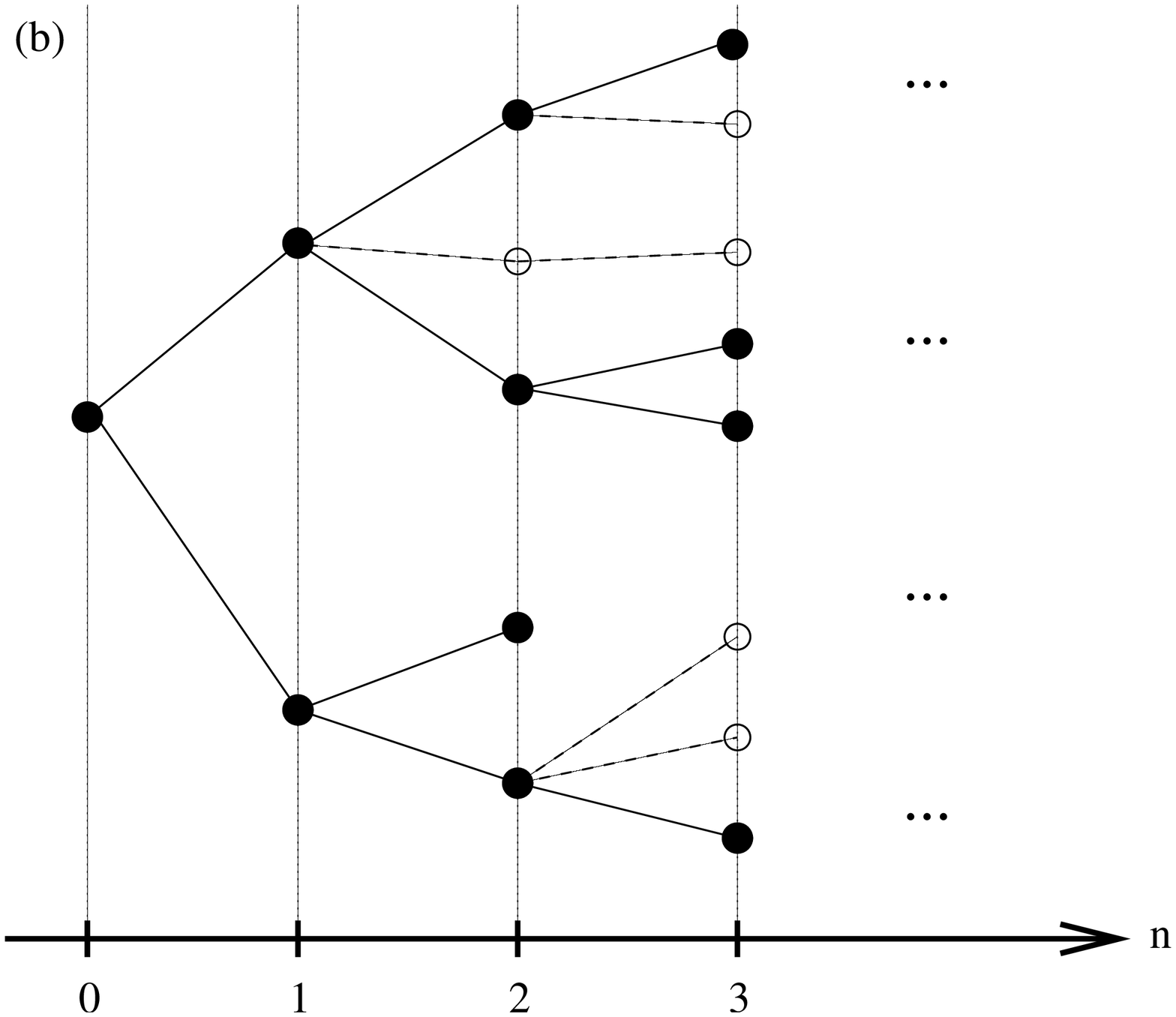}

\caption{Probabilistic realization of branching process $\left\{ Z_{n}\right\} _{n=0}^{\infty}$
without vaccination (a) and with vaccination (b)}
\end{figure}

Let $\beta_{n}$ be the proportion of the infected individuals who
infect $n$ individuals. Hence, $\beta_{0}$ is the proportion of
infected individuals who do not infect anyone before they recover
or die. If $\beta_{0}=0$, then $Z_{n}>Z_{0}$ for all $n$ and the
process is not interesting, so we assume that $\beta_{0}>0$. Likewise
if $\beta_{0}=1$ then $\pi_{n}$, the probability of extinction of
the virus after $n$ generations, will be clearly equal to 0 for all
$n$. Hence, we assume $0<\beta_{0}<1$ throughout. Clearly, $\pi_{n}<\pi_{n+1}$
for all $n$. If the infection is extinct in generation $n$, it will
remain extinct for $n+1$ and all following generations. Then, $\left\{ \pi_{n}\right\} _{n=1,2,...}$
is a monotone increasing sequence, so we let $\pi=\lim_{n\rightarrow\infty}\pi_{n}$.
Here $\pi$ is the probability that the infection will eventually
die out in the population. Suppose $F(x)=\beta_{0}+\beta_{1}x+\beta_{2}x^{2}+...$
then $F(\pi_{n})=\beta_{0}+\beta_{1}\pi_{n}+\beta_{2}\pi_{n}^{2}+...$
$=\sum\beta_{r}\pi_{n}^{r}$. $F(\pi_{n})$ is the sum of the probabilities
of the $r+1$ distinct ways in which the infection can become extinct. 

In the next section we introduce probability generating function $I(x)$
for $\{\theta_{i}\}$ and prove results similar to the one discussed
above on $F(x)$ in broader way. We study uniform convergence properties
and related consequences of $I(x)$ for $x\in(-\mathbb{R},\mathbb{R})$
and $x\in(-1,1).$ The results on uniform convergence of $I(x)$ provided
in this work are not new because $I(x)$ is a power series (introduced
in the next section). We study extinction probabilities \emph{with
respect to} the \emph{}mean number of infectious individuals generated
by starting $Z_{0}=1$. This mean number is defined as $I'(1).$ In
the mathematical epidemiology literature, $I'(1)$ would be referred
to as the basic reproductive number. The process $\{ Z_{n}\}$ is
called supercritical if $I'(1)>1$ and the process is called sub critical
if $I'(1)\leq1.$ We have not seen rigorous mathematical treatment
of the branching process for the spread of virus in the population
(except elementary branching process applications for understanding
measles outbreaks \cite{key-0} and for virus HIV spread within the
host cell dynamics\cite{key-3}). Extinction probability results in
the pre-vaccination (i.e. without vaccination) scenario are derived
from branching process studies available on offspring distribution,
surnames and other general branching processes analysis (see for example,
\cite{key-1,key-2,key-4,key-5,key-7,key-8,key-9,key-10,key-11} and
fundamental results on uniform convergence can be seen in many classical
analysis books). Hence some results in section 2 do not bring new
mathematics, but rather new applications. In section 3, we bring novelty
through theoretical arguments on EBP when vaccination is introduced
in the population. By vaccinating the individuals, the resultant probability
density functions of infection densities shifts to the left and extinction
of virus is quicker. We also explore some consequences of vaccination.

\section{Elementary results on extinction probability}

We saw $\Pr\left[N_{r}=i\right]=\theta_{i}$ and $N_{r}$ are independent
discrete random variables. Suppose $I_{2}$ and $J_{r}$ are probability
generating functions of $Z_{2}$ and $N_{r}$ for $r=1,2,3\cdots.$
Then we know $J_{1}=I_{1}$ because $Z_{1}=N_{1},$ $J_{1}.J_{2}=I_{2}.$
Suppose distribution of $N_{r}$ are same (say some $J$), then $I_{2}=\left[J\right]^{2}$
where $I_{2}=I_{2}(x)=\sum_{i=0}^{\infty}\Pr\left[Z_{2}=i\right]x^{i}$
is the probability generating function of $Z_{2}.$ In general let
$I_{n}$ is probability generating function of $Z_{n}.$ Since $Z_{2}=N_{1}+N_{2}+\cdots+N_{Z_{1}}$,
it follows that $\Pr\left[Z_{2}=k+1/Z_{1}=k\right]=\Pr\left[\sum_{i=1}^{k}N_{i}=k+1\right].$
We can prove $I_{2}=I_{1}(I_{1})$ and $I_{n}=I_{n-1}(I_{1})$ \cite{key-11}.
Let $I(x)=\sum_{i=0}^{\infty}\Pr\left[N=i\right]x^{i}$ $(\left|x\right|<1)$
be the probability generating function for $\{\theta_{i}\}.$ We prove
convergence properties of $I(x)$ and hence bring elementary results
on extinction of probability of the virus from the population. 

\begin{thm}
\label{theorem-1}Suppose $I(x)$ converges for $\left|x\right|<\mathbb{R}$,
then for a $\delta>0,$ $I(x)$ converges uniformly on $\left[-\mathbb{R}+\delta,\mathbb{R}-\delta\right]$
and $I(x)$ is continuous and infinitely differentiable on any open
interval $(-\mathbb{R},\mathbb{R}).$
\end{thm}
\begin{proof}
It is enough to show existence of $I'$ i.e. differentiability of
$I$, then it follows continuity of $I$. This is because,

\begin{eqnarray*}
lim_{x\rightarrow p}\left[I(x)-I(p)\right] & = & lim_{x\rightarrow p}\frac{I(x)-I(p)}{x-p}\left[x-p\right]=I'(x).0=0\end{eqnarray*}

Hence $I$ is continuous. $I'(x)=\sum i\theta_{i}x^{i-1}$ ($\left|x\right|<\mathbb{R}$).
Also, as $i\rightarrow\infty$ , it follows that $\limsup\sqrt[i]{\left|\theta_{i}\right|}$$=\limsup\sqrt[i]{i\left|\theta_{i}\right|}.$
This is due to fact that $i\rightarrow\infty$ , it follows that $\sqrt[i]{i}=1$.
Since $\limsup\sqrt[i]{\left|\theta_{i}\right|}$ as $n\rightarrow\infty$,
it follows that $\sum\theta_{i}$ is convergent. Therefore, $\sum\theta_{i}x^{i}$
and $\sum i\theta_{i}x^{i}$ have interval of convergence in $(-\mathbb{R},\mathbb{R}).$
Since $I(x)$ is a power series, the required results is straight
forward and can be seen in several classical analysis books. Suppose
the series $\sum_{n=0}^{\infty}a_{n}x^{n}$ converges for $\left|x\right|<\mathbb{R},$
and define $g(x)=\sum_{n=0}^{\infty}a_{n}x^{n}$ $(\left|x\right|<\mathbb{R}).$
Then $g$ converges uniformly on $[-\mathbb{R}+\epsilon,\mathbb{R}-\epsilon],$
for every $\epsilon>0.$ Here $g$ is continuous and differentiable
in $(-\mathbb{R},\mathbb{R})$ and $g'(x)=\sum_{n=1}^{\infty}na_{n}x^{n-1}.$ 

The above theorem established that the probability generating function
of the infected distribution is differentiable over the real line
if $\left|x\right|<\mathbb{R}$ holds. However conventionally, we
define $I(x)$ for $\left|x\right|<1.$
\end{proof}
\begin{thm}
For a given $\epsilon,\delta>0$, let $\mathcal{B}_{\epsilon}(\sum\theta_{i})=\left\{ x\in\mathbb{R}:\left|x-\sum\theta_{i}\right|<\epsilon\right\} $
and $\mathcal{B}_{\delta}(1)=$ $\{ x\in\mathbb{R}:\left|x-1\right|$
$<\delta\}$ and $\theta_{i}<\theta_{i+1}$ for $\{ i=0,1,2\cdots\}.$
Then for every $\mathcal{B}_{\epsilon}(\sum\theta_{i})$, there exists
a $\mathcal{B}_{\delta}(1)$ with the property that for all $x\in\mathcal{B}_{\delta}(1)$,
$x\neq1$, it follows that $I(x)\in\mathcal{B}_{\epsilon}(\sum\theta_{i}).$
\end{thm}
\begin{proof}
Let $T_{i}=\sum\theta_{i}.$ Consider the sequence $\left\{ T_{i}\right\} _{i=0,1,2,}$.
Here $T_{0}=\theta_{0},$ $T_{1}=\theta_{0}+\theta_{1},$ $T_{2}=\theta_{0}+\theta_{1}+\theta_{2}$
and so on. We assume $\theta_{0}>0.$ Otherwise, if $\theta_{0}=0$
the branching process will terminate. Since $\theta_{0}>0$, and from
the hypothesis we have, $\theta_{1}>\theta_{0}>0$ and $\theta_{2}>\theta_{1}$
and so on. Also, we know $\theta_{i}<1$ because $\{\theta_{i}\}$
forms a complete probability distribution. Therefore $\{ T_{i}\}$
is convergent. We have the very important result by Abel, which states
that when $\sum\theta_{n}$ is convergent and as $x\rightarrow1$,
it follows that $\lim\sum\theta_{n}x^{n}=\sum\theta_{n}$ $(-1<x<1)$.
Hence from this result it follows that $x\in\mathcal{B}_{\delta}(1)$$\Rightarrow I(x)\in\mathcal{B}_{\epsilon}(T_{i}).$
\end{proof}
\begin{rem}
The elementary theorem \ref{ele-th3} and the arguments in the proof
given can also be found for studying extinction of surnames of families,
extinction of particular offspring, particle extinction in a volume
of gas and other biological applications. Otherwise, this theorem
is very well known, but we state it here in the language of infectious
disease extinction in the population. This will help readers to compare
the above elementary result with the new results in the next section
when a vaccine that prevents disease transmission is introduced into
the population. 
\end{rem}
\begin{thm}
\label{ele-th3}Given $I(x)$ for $x\in[0,1]$ and size of infected
population in zeroth generation is one. If the process is sub critical
or critical then $\lim\Pr\left[Z_{n}=0\right]=1$ as $n\rightarrow\infty$
and if the process is supercritical then there exists a $\pi(\neq1)\in[0,1]$
such that $I(\pi)=\pi$ and as $n\rightarrow\infty$ it follows that
$\lim\Pr\left[Z_{n}=0\right]=\pi.$
\end{thm}
\begin{proof}
Consider the closed interval $[0,1].$ We saw in theorem \ref{theorem-1}
that $I(x)$ is continuous in $[0,1]$ and also $I'(x)=\sum i\theta_{i}x^{i}>0.$
Therefore, $I'(x)$ is strictly increasing in $[0,1]$. If the system
is sub critical then $I'(1)\leq1.$ This means $I'(x)<1$ for $x\in[0,1).$
This implies, $\int_{x}^{1}I'(x)dx<\int_{x}^{1}dx$ $\Rightarrow$
$1-I(x)<1-x\Rightarrow I(x)>x.$ $I(1)=\sum\theta_{i}=1.$ This means
the curve of $I(x)$ never touches $I(x)=x$ in $[0,1)$ and at $1$
it will be $1.$ Therefore, $I(1)=1$ has to be a unique fixed point
for $I$ in $[0,1].$Now it is left for us to show that $\lim\Pr\left[Z_{n}=0\right]=1.$
We can easily prove that $\left\{ \Pr\left[Z_{n}=0\right]\right\} $
is convergent. Let $\pi$ be the limit of this sequence, so that as
$n\rightarrow\infty,$ it follows that $\lim\Pr\left[Z_{n}=0\right]=\pi.$
Also, $\lim I\left(\Pr\left[Z_{n-1}=0\right]\right)=I(\pi)$ and $I(\pi)=\pi.$
But we saw above that, when $I'(1)\leq1,$ then $I(1)=1$ is unique
fixed point. Hence $\lim\Pr\left[Z_{n}=0\right]=1.$

Since $I'(x)=\sum i\theta_{i}x^{i-1}$>0 for $x\in[0,1],$ $I'(x)$
is strictly increasing, then there exists a $a\in(0,1)$ such that
for $b\in(a,1)$ we will have 

\begin{eqnarray}
1 & < & I'(b)<I'(1)\label{two}\\
\Rightarrow\int_{b}^{1}dx & < & \int_{b}^{1}I'(x)dx\nonumber \\
\Rightarrow1-b & < & 1-I'(b)\Rightarrow b>I(b)\quad\textrm{for $b\in(a,1)$}\label{three}\end{eqnarray}
These kind of above arguments can be found in basic setting of branching
process (see \cite[(chapter 12)]{key-1}, \cite[(chapter 4)]{key-5},
\cite[(chapter 1)]{key-10}, \cite[(chapter 6)]{key-11}, \cite[(chapter 0)]{key-15}).
Equation \ref{two} is a situation of supercritical process. When
process obeys such property then our aim is to show that there exists
two fixed points $\pi$ and $1$ in $[0,1]$ such that $\pi\neq1.$
Every value of $b$ between $a$ and $1$ (exclusive), the value of
line $y=x$ is greater than the value of the curve $y=I(x)$ for $x\in(a,1)$
on the $XY-$ plane. Therefore, this situation is not conducive for
us to find a $\pi$ such that $I(\pi)=\pi.$ If we bring a relation
$I(x)-x=0$ for some $x<1$ then prove that $x$ is unique, then we
are ready to compute extinction probabilities. Consider $I(x)-x.$
We have $I(0)-0=\theta_{0}>0$ and from equation \ref{three} $I(b)<b$
for $b\in(a,1).$ This means $I(x)-x<0$ for $x=b.$ Thus, 

\begin{eqnarray*}
I(x)-x & = & \left\{ \begin{array}{c}
>0\quad\textrm{for}\: x=0\\
<0\quad\textrm{for}\: x=b.\end{array}\right.\end{eqnarray*}

The intermediate value theorem says if any continuous function on
a given closed interval with values ranging from negative to positive,
then this function will be zero for some value in between the same
closed interval in the domain. Therefore $I(x)-x=0$ for some $x\in(0,b).$
Let $I(\pi)=\pi$ for $\pi\in(0,b).$ We are sure that in $(a,1)$
there will be no fixed point. Suppose $w$ is the another fixed point
in $(0,1)$ then either $w\in(0,\pi)$ or $w\in(\pi,1).$ Since $w$
is another fixed point, $I(w)-w=0,$ $\pi$ is also a fixed point,
so $I(\pi)-\pi=0$ and $I(1)-1=0.$ $I(x)-x$ is differentiable and
continuous in $(0,1),$ so

\begin{eqnarray}
\left\{ \begin{array}{ccc}
I(w_{1})-w_{1} & = & 0\quad\textrm{for}\:0<w_{1}<w\\
I(w_{2})-w_{2} & = & 0\quad\textrm{for}\: w<w_{2}<\pi\end{array}\right.\label{four}\end{eqnarray}

and 

\begin{eqnarray}
\left\{ \begin{array}{ccc}
I(w_{1})-w_{1} & = & 0\quad\textrm{for}\:\pi<w_{1}<w\\
I(w_{2})-w_{2} & = & 0\quad\textrm{for}\: w<w_{2}<1\end{array}\right.\label{five}\end{eqnarray}

These equations \ref{four} \& \ref{four} together contradict the
fact that $I(x)$ in strictly increasing in $(0,1).$ Therefore $w=\pi.$
This means $\lim\Pr\left[Z_{n}=0\right]=\pi.$ Hence we have two fixed
points $\pi$ and $1$ in $x\in[0,1].$
\end{proof}
\begin{rem}
If the mean number of infections exceed unity, then the virus from
the population will go extinct in a finite number of generations with
probability $\pi.$ Suppose $Z_{0}=N$ then each of these individual's
process can be treated independently. Thus $\lim\Pr\left[Z_{n}=0\right]=\pi^{N}$
and $\lim\Pr\left[Z_{n}=\infty\right]=1-\pi^{N}.$ $\lim\Pr\left[Z_{n}=\infty\right]\rightarrow0$
when $N\rightarrow\infty$even if $\pi>>0.$ This means virus will
go extinct even if population at the zeroth stage has a large number
of infected individuals. 
\end{rem}
\begin{prop}
Let $\{ Z_{n,\phi}\}$ is the size of the infected population during
post vaccine scenario. If $\{ Z_{n}\}$ is convergent then $\{ Z_{n,\phi}\}$
is also convergent, but converse is not true for $\phi>0.$ 
\end{prop}
\begin{proof}
Suppose $\{ Z_{n}\}$ is convergent. Then for $\delta>0$, $\left|Z_{m}-Z_{n}\right|<\delta$
whenever $m,n\geq N\in\mathbb{N}.$ We have,

\begin{eqnarray*}
Z_{n+1,\phi} & = & N_{1,\phi}+N_{2,\phi}+\cdots+N_{Z_{n},\phi}\\
 & \leq & N_{1}+N_{2}+\cdots+N_{Z_{n}}=Z_{n+1}\end{eqnarray*}

Therefore $\left|Z_{m,\phi}-Z_{n,\phi}\right|<\delta$ whenever $m,n\geq N\in\mathbb{N}.$
Hence $\{ Z_{n,\phi}\}$ is convergent. 
\end{proof}

\section{Spreading is restricted by vaccination}

Suppose the process $\{ Z_{n}\}$ is controlled by introducing vaccination
into the population, which protects susceptible individuals from infection
and therefore reduces the number of secondary infections produced
per infected individual. We model the probability density function
of the number of secondary infections per infected person as a gamma
function, so that it covers a wide range of epidemiologically-plausible
scenarios. Vaccination will shift the peak of the probability density
function to the left. Our aim here is to estimate the corresponding
change in the time to extinction brought about by vaccination. 

\begin{prop}
\label{proposiotion 3.1} Suppose the time to extinction without vaccination
is $\tau_{a}$ and $P[N=w_{a}]=\Gamma_{a}$ $[w_{a},T_{a},\alpha_{a}]$
and suppose the time to extinction with vaccination is $\tau_{b}$
and $P[N=w_{b}]=\Gamma_{b}[w_{b},T_{b},\alpha_{b}]$. If mean of $\Gamma_{b}\textrm{ }$
is less than the mean of $\Gamma_{a}$ then $\tau_{b}<\tau_{a}.$
\end{prop}
\begin{proof}
Let $\overline{w_{a}}$, $\overline{w_{b}}$ means of $\Gamma_{a}$
and $\Gamma_{b}$. We know that by vaccinating the population, the
new number of infections generated by one infected is reduced and
hence the sizes of $\{ Z_{n}\}$ at each stage $n$ is affected. We
have $\overline{w_{b}}<\overline{w_{b}}\Rightarrow T_{b}\alpha_{b}<T_{a}\alpha_{a}.$
Hence the time to extinct will be earlier with vaccine than that of
without vaccination. 
\end{proof}
We discuss here some further consequences of vaccination. Suppose
$\sum_{w_{a}=0}^{\infty}w_{a}\theta_{w_{a}}=$ $\int_{-\infty}^{\infty}w_{a}$
$d(\Gamma_{a}$ $[w_{a},T_{a},\alpha_{a}])$ $=T_{a}\alpha_{a}.$
Here $d(\Gamma[.])$ is gamma distribution function. We have seen
in the previous section that if $\sum w_{a}\theta_{w_{a}}\leq1$ then
the chance of extinction is one and if $\sum w_{a}\theta_{w_{a}}>1$
then the chance is  $\pi$ for $\pi<1.$ This also means $\lim\Pr\left[Z_{n_{a}}=0\right]=1$
as $n_{a}\rightarrow\infty$. Suppose $\sum_{w_{b}=0}^{\infty}w_{b}\theta_{w_{b}}=\int_{-\infty}^{\infty}w_{b}d(\Gamma_{b}[w_{b},T_{b},\alpha_{b}])=T_{b}\alpha_{b}.$
Please note $w_{a}$ and $w_{b}$ denote infected individuals and
$\theta_{w_{a}}$ and $\theta_{w_{b}}$ denote infected population
distributions for pre and post vaccination. $\{ w_{j}\},\, j=1,2,3,...$
are individuals in post vaccination scenario. Here $\lim\Pr\left[Z_{n_{b}}=0\right]=1$
as $n_{b}\rightarrow\infty.$ Note that $Z_{n_{b}}<Z_{n_{a}}$ and
$n_{b}<n_{a}$ (by proposition \ref{proposiotion 3.1}) If we consider
probability densities of these two scenarios on the $XY-$plane then
the mean of the infection distribution with vaccination shifts to
left of the mean of the infection distribution without vaccination
i.e. $T_{b}\alpha_{b}<T_{a}\alpha_{a}.$ Suppose $\{ T_{j}\},\{\alpha_{j}\}$
for $j=1,2,...,n$ are set of parameters in the plausible range of
post vaccine scenario. If $\mathbf{T}=Min\{ T_{1},...,T_{n}\}$ and
$\mathbf{\alpha}=Min\{\alpha_{1},...,\alpha_{n}\}$ then post vaccination
infection density will have one of these three parameters combinations:
$I.\,\{\mathbf{T},\alpha_{j}\}$ or $II.\,\{ T_{j},\mathbf{\alpha}\}$
or $III.\,\{\mathbf{T},\mathbf{\alpha}\}.$ In the cases of $I$ and
$II$ the infection density could be realistic, but $III$ is not
certainly an expected density. To avoid situation $-III$, we choose
a set $E_{j}$ consisting of all possible combinations of parameters
$\{ T_{j},\alpha_{j}\}$ such that $\Gamma-$ density is never a decreasing
function. 

The corresponding inequalities due to situations $I$, $II$ and $III$
are $\mathbf{T}\alpha_{j}<T_{a}\alpha_{a},$ $T_{j}\mathbf{\alpha}<T_{a}\alpha_{a}$
and $\mathbf{T\alpha}<T_{a}\alpha_{a}$ for some $j\textrm{ }$ such
that post vaccination parameter sets are in $E_{j}.$ Since $\overline{w_{b}}<\overline{w_{a}}$
it can be viewed that $\tau_{b}<\tau_{a}.$ However the mean of $\Gamma_{j}$
for any $j$ is less than the mean of $\Gamma_{a}.$ This implies,

\begin{eqnarray*}
\Pr[Z_{n_{a}}=0](\textrm{as }n_{a}\rightarrow\infty) & < & \Pr[Z_{n_{b}}=0](\textrm{as }n_{b}\rightarrow\infty)\end{eqnarray*}

Also, $\Pr[Z_{n_{a}}=0]<[Z_{n_{j}}=0].$ Note that $n_{a},n_{j}$
are not associated with the same time axis because pre and post vaccination
situations cannot occur at the same time. In section 2, $Z_{n},Z_{n+1},...$
are associated to the same time axis. Consider the sequence $\left\{ \Pr\left[Z_{n_{a}}=0\right]\right\} _{n_{a}=0,1,2\cdots}.$
Note that the suffix $a$ is used here to indicate that size of the
population in pre vaccine scenario. We have $I(0)=\Pr\left[Z_{1}=0\right]=\theta_{0}>0$
and also $\theta_{0}<1$. Then by the mathematical induction and monotonic
property of $I$ we can show that \begin{eqnarray*}
I(\Pr\left[Z_{(n-1)_{a}}=0\right]) & < & I(\Pr\left[Z_{n_{a}}=0\right])<I(1)\\
\Rightarrow\Pr\left[Z_{n_{a}}=0\right] & < & \Pr\left[Z_{(n+1)_{a}}=0\right]<1\end{eqnarray*}

Similarly by considering the sequence $\left\{ \Pr\left[Z_{n_{j}}=0\right]\right\} _{n_{j}=0,1,2\cdots}$
we can show that

\begin{eqnarray*}
I(\Pr\left[Z_{(n-1)_{j}}=0\right]) & < & I(\Pr\left[Z_{n_{j}}=0\right])<I(1)\\
\Rightarrow\Pr\left[Z_{n_{j}}=0\right] & < & \Pr\left[Z_{(n+1)_{j}}=0\right]<1\end{eqnarray*}

Thus $\left\{ \Pr\left[Z_{n_{a}}=0\right]\right\} $ and $\left\{ \Pr\left[Z_{n_{j}}=0\right]\right\} $
are convergent sequences. Moreover $\left\{ \Pr\left[Z_{n_{j}}=0\right]\right\} $
converges earlier due to vaccination. 

\begin{thm}
Let $\theta_{w_{j}}$ be the probability that an individual infects
$w_{j}$ individuals during post vaccination, then $\{\sum_{j}(\sum_{w_{j}}w_{j}\theta_{w_{j}})\}^{2}\leq(\sum_{j}T_{j}^{2})(\sum_{j}\alpha_{j}^{2}).$
\end{thm}
\begin{proof}
Since $\{ T_{1},...,T_{n}\}$ and $\{\alpha_{1},...,\alpha_{n}\}$
are non-zero real numbers, we have from the Cauchy-Schwartz inequality,

\begin{eqnarray*}
(\sum_{j}T_{j}^{2})(\sum_{j}\alpha_{j}^{2}) & \geq & (\sum_{j}T_{j}\alpha_{j})^{2}\\
 & = & \sum_{j}\{\int_{-\infty}^{\infty}w_{j}d(\Gamma_{j}[w_{j},T_{j},\alpha_{j}])\}^{2}\\
 & = & \sum_{j}\{(\sum_{w_{j}}w_{j}\Pr[N=w_{j}])\}^{2}\\
 & = & \{\sum_{j}(\sum_{w_{j}}w_{j}\theta_{w_{j}})\}^{2}\end{eqnarray*}

Thus this inequality realtes the mean number of infections and the
parameter values in post vaccination. Hence for modeling super critical
and sub critical situations with respect to parameter values this
inequality can be utilised. 
\end{proof}
\begin{prop}
Suppose the mean distance between the peaks of the infection densities
is $\overline{d}(P_{j},P_{a}).$ Here $P_{a}(w_{a},\Gamma_{a}(w_{a}))$
are $P_{j}(w_{j},\Gamma_{j}(w_{j}))$ are the peak points for the
pre and post vaccination densities. Then $\Gamma(\overline{w_{j}})=\Gamma(w_{a})\pm\left[\frac{1}{n^{2}}\left\{ \sum_{j}d(P_{j},P_{a})\right\} ^{2}-(w_{a}-\overline{w_{j}})^{2}\right]^{\frac{1}{2}}.$ 
\end{prop}
\begin{proof}
Let $P_{j}(w_{j},\Gamma_{j}(w_{j}))$ be the point representing the
peak of the density function after vaccination, in $E_{j}.$ If $P_{a}(w_{a},\Gamma_{a}(w_{a}))$
is the point representing the peak of the density function before
vaccination, then Euclidean distance between these two points is $d(P_{j},P_{a})$
(see Figure \ref{figure-peaks} ). The mean ( $\overline{d}(P_{j},P_{a})$)
of this function over all possible $j$ is 

\begin{eqnarray*}
\overline{d}(P_{j},P_{a}) & = & \frac{\sum_{j}d(P_{j},P_{a})}{n}\end{eqnarray*}

\begin{figure}
\includegraphics[%
  scale=0.5]{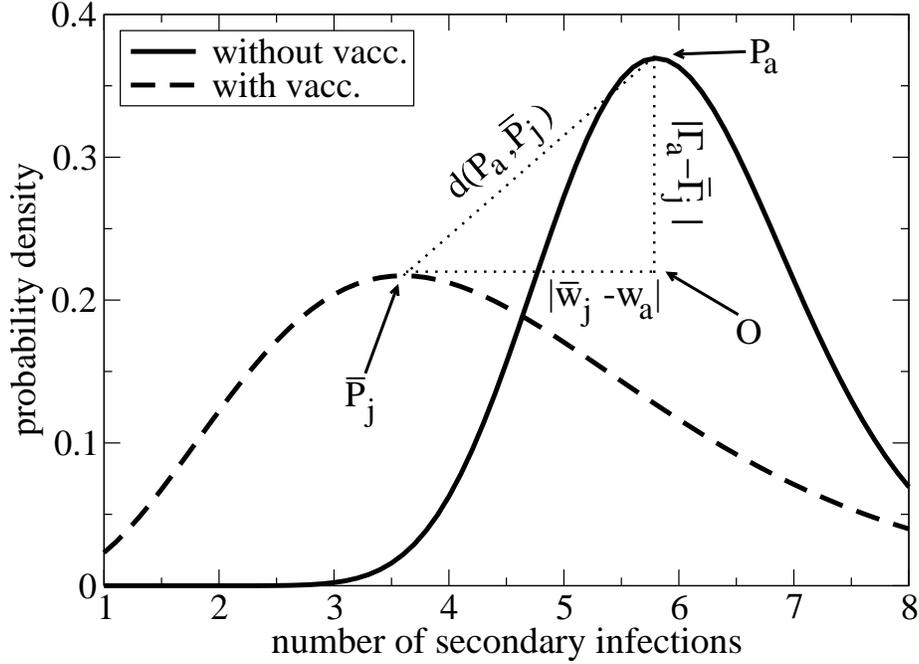}

\caption{\label{figure-peaks}Examples of infection peaks with ($P_{j}$)
and without ($P_{a}$) vaccination. For the with-vaccination scenario,
the average peak over many possible types is shown (see text). }
\end{figure}

For $j$ possible types of post vaccine densities, let $\overline{P_{j}}(\overline{w_{j}},\Gamma(\overline{w_{j}}))$
be the point corresponding to the mean peak. $O$ is a point which
is at a distance of $w_{a}-\overline{w_{j}}$ from $\overline{P_{j}}$
to the right. Suppose $\overline{d}(P_{a},\overline{P_{j}})$ is the
distance from $P_{a}$ to $\overline{P_{j}}$. Then we know from Pythogoren
right triangle principle that $d(P_{a},\overline{P_{j}})=\sqrt{(OP_{a})^{2}+(O\overline{P_{j}})^{2}}.$
From general principles of means, we can easily verify the fact that
$d(P_{a},\overline{P_{j}})=\overline{d}(P_{j},P_{a}).$ Therefore,

\begin{eqnarray*}
\frac{1}{n^{2}}\left\{ \sum_{j}d(P_{j},P_{a})\right\} ^{2} & = & (w_{a}-\overline{w_{j}})^{2}+(\Gamma(w_{a})-\Gamma(\overline{w_{j}}))^{2}\\
\Gamma(\overline{w_{j}}) & = & \Gamma(w_{a})\pm\left[\frac{1}{n^{2}}\left\{ \sum_{j}d(P_{j},P_{a})\right\} ^{2}-(w_{a}-\overline{w_{j}})^{2}\right]^{\frac{1}{2}}\end{eqnarray*}

This proposition relates the shape of the distribution of the number
of secondary infections per infected individual before vaccination
to the same after vaccination. This may be useful in situations where
both distributions have small skew and $w_{a}$, $w_{j}$ correspond
to the expectation values of the distributions.
\end{proof}
\begin{cor}
Suppose $I_{a}(\sigma^{2})$ and $I_{j}(\sigma^{2})$ are variances
of infections corresponding to the distributions $\{\theta_{n_{a}}\}$
and $\{_{\theta_{n_{j}}}\}.$ Then we have $I_{a}(\sigma^{2})-I_{j}(\sigma^{2})=$

\begin{eqnarray*}
 &  & \sum_{n_{a}}\left(n_{a}-I_{a}'(1)\right)^{2}\Pr\left[N=n_{a}\right]-\sum_{n_{j}}\left(n_{j}-I_{j}'(1)\right)^{2}\Pr\left[N=n_{j}\right]\\
 & = & \sum_{n_{a}}n_{a}^{2}\Pr\left[N=n_{a}\right]+2\left\{ \left(I_{j}'(1)\right)^{2}-\left(I_{a}'(1)\right)^{2}\right\} -\sum_{n_{j}}n_{j}^{2}\Pr\left[N=n_{j}\right]+\\
 &  & \qquad\qquad\qquad\qquad\qquad\qquad\qquad\qquad\left\{ \left(I_{a}'(1)\right)^{2}-\left(I_{j}'(1)\right)^{2}\right\} \\
\\ & = & \sum_{n_{a}}n_{a}^{2}\Pr\left[N=n_{a}\right]-\sum_{n_{j}}n_{j}^{2}\Pr\left[N=n_{j}\right]-2\left\{ \left(I_{a}'(1)\right)+\left(I_{j}'(1)\right)\right\} \overline{d}_{\left|A\right|}+\\
 &  & \qquad\qquad\qquad\qquad\qquad\qquad\qquad\qquad\left\{ \left(I_{a}'(1)\right)+\left(I_{j}'(1)\right)\right\} \overline{d}_{\left|A\right|}\\
\\ & = & \sum_{n_{a}}n_{a}^{2}\Pr\left[N=n_{a}\right]-\sum_{n_{j}}n_{j}^{2}\Pr\left[N=n_{j}\right]-\left\{ \left(I_{a}'(1)\right)+\left(I_{j}'(1)\right)\right\} \overline{d}_{\left|A\right|}\end{eqnarray*}

\end{cor}
In the above $I_{a}'(1)-I_{j}'(1)=\overline{d}_{\left|A\right|}.$
The above relation indicates that difference between variability between
infection process also depends upon the adjusted distance between
the peaks described above. If the process is critical then the above
difference is expressed as a straight combination of second central
moments. 

\begin{rem}
\begin{eqnarray*}
\sum n_{a}\Pr[Z_{n_{a}}=0]\leq1 & \Rightarrow & \sum n_{j}\Pr[Z_{n_{j}}=0]\leq1\\
\Rightarrow\overline{d}(P_{j},P_{a}) & \leq & v^{2}\;\textrm{for some $v\leq1$}\end{eqnarray*}

Again,

\begin{eqnarray*}
\sum n_{a}\Pr[Z_{n_{a}}=0]>1 & \Rightarrow & \left\{ \begin{array}{c}
\sum n_{j}\Pr[Z_{n_{j}}=0]>1\;(\textrm{no impact of vaccination})\\
\sum n_{a}\Pr[Z_{n_{a}}=0]\leq1\;(\textrm{impact of vaccination})\end{array}\right.\end{eqnarray*}

Since vaccine has a positive role in reducing the infection process
such that as $n_{j}>N$ (for some large $N$) it implies $\sum n_{j}\Pr[Z_{n_{j}}=0]\leq1.$
When the process $\{ Z_{n_{a}}\}$ is supercritical then initially
the process $\{ Z_{n_{j}}\}$ could also be supercritical and eventually
$\{ Z_{n_{j}}\}$ will attain sub criticality or criticality. Hence
even if $\sum n_{a}\Pr[Z_{n_{a}}=0]>1$, eventually $d>1.$ The larger
the value of $d$ the larger is the impact of vaccination and extinction
occures much earlier with probability $\pi_{j}$$(\pi_{j}<\pi).$
\end{rem}

\section{Conclusions}

Branching process analyses have provided a rigorous perspective on
numerous physical problems in the past. Here, we have introduced the
foundations for application of branching processes to studying the
early stages of epidemic spread in a population, both with and without
vaccination. We modified the standard branching process to account
for the effects of vaccination, showing how vaccination will decrease
the time to extinction. There remains much work to be done in epidemic
branching processes, particularly in the incorporation of other types
of control measures (quarantine) and in the incorporation of realistic
population structures such as age-structure.

\end{document}